\documentclass[11pt]{amsart}
\usepackage{amsmath}
\usepackage{amssymb}
\usepackage{amsfonts}
\usepackage{graphicx}

\DeclareSymbolFont{AMSb}{U}{msb}{m}{n}
\DeclareMathSymbol{\N}{\mathbin}{AMSb}{"4E}
\DeclareMathSymbol{\Z}{\mathbin}{AMSb}{"5A}
\DeclareMathSymbol{\R}{\mathbin}{AMSb}{"52}
\DeclareMathSymbol{\Q}{\mathbin}{AMSb}{"51}
\DeclareMathSymbol{\I}{\mathbin}{AMSb}{"49}
\DeclareMathSymbol{\C}{\mathbin}{AMSb}{"43}

\theoremstyle{definition}

\theoremstyle{corollary}

\theoremstyle{example}

\theoremstyle{note}

\theoremstyle{notation}

\numberwithin{equation}{section}
\begin{document}
\title[COMBINATORIALLY RICH SETS NEAR ZERO]
{DYNAMICAL CHARACTERIZATIONS OF COMBINATORIALLY RICH SETS NEAR ZERO}
\author{Sourav Kanti Patra}
\address{Sourav Kanti Patra, Department of Mathematics, Ramakrishna Mission Vidyamandira,
Belur Math, Howrah-711202, West Bengal, India}
\email{souravkantipatra@gmail.com}
%\author{}
%\address{}
%\email{}
%\thanks{}
\keywords{Algebra in the Stone-$\breve{C}$ech compactification,
Dynamical system, Uniform recurrence near zero, Proximality near zero,
JIUR, JIAUR, JIUR near zero, JIAUR near zero}

\begin {abstract}
Hindman and Leader first introduced the notion of Central sets 
near zero for dense subsemigroups of $((0,\infty),+)$ and proved 
a powerful combinatorial theorem about such sets. Using the 
algebraic structure of the Stone-$\breve{C}$ech compactification, 
Bayatmanesh and Tootkabani generalized and extended this combinatorial 
theorem to the central theorem near zero. Algebraically one can define 
quasi-central set near zero for dense subsemigroup of $((0,\infty),+)$, 
and they also satisfy the conclusion of central sets theorem near zero. 
In a dense subsemigroup of $((0,\infty),+)$, C-sets near zero are the 
sets, which satisfies the conclusions of the central sets theorem 
near zero. Like discrete case, we shall produce dynamical characterizations 
of these combinatorically rich sets near zero.

AMS subjclass [2010] : 37B20; 37B05; 05B10.
\end{abstract}

\maketitle

\section {introduction}
Furstenberg, defined the concept of a central subset of positive 
integers [4, Definition 8.3] and proved several important properties 
of such sets using notions from topological dynamics.\\

\textbf{Definition 1.1} A dynamical system is a pair $(X,\langle T_s\rangle _{s \in S })$ 
such that,

(i) $X$ is compact Hausdorff space,

(ii) $S$ is a semigroup,

(iii) for each $s\in S$, $T_s:X\rightarrow X$ and $T_s$ is continuous, and

(iv) for all $s,t$ , $T_s\circ T_t=T_{st}$.\\

Inspired by the fruitful interaction between Ramsey theory and ultrafilters 
on semigroups, Bergelson and Hindman, with the assistance of B. Weiss, later 
proved on algebraic characterization of central sets in 
$\mathbb{N}$ [3, Section 6]. Using this algebraic characterization a as a 
definition enabled them easily to extend the notion of a central set to any semigroup.\\
Let us now give a brief description of the algebraic structure of $\beta S$ for a 
discrete semigroup $(S,\cdot)$. We take the points of $\beta S$ to be the ultrafilters 
on $S$, identifying the principal ultrafilters with the points of $S$ and thus 
pretending that $S\subseteq \beta S$. Given $A\subseteq S$ let us 
set, $\bar{A}=\{ p\in \beta S : A\in p\}$. Then the set $\{ \bar{A} : A\subseteq S\}$ 
is a basis for a topology on $\beta S$. The operation $\cdot$ on $S$ can be extended to the 
Stone-cech compactification $\beta S$ of $S$ so that $(\beta S,\cdot)$ is a compact right 
topological semigroup(meaning that for any $p\in \beta S$), the function $\rho_p:\beta S$ 
defined by, $\rho_p(q)=q\cdot p$ is continuous) with $S$ contained in its topological center 
(meaning that for any $x\in S$ the function $\lambda_x : \beta S \rightarrow \beta S$ 
defined by $\lambda_x (q)=x\cdot q$ is continuous). Given $p,q \in \beta S $ and 
$A \in S$ $A \in p\cdot q$ if and only if $\{ x \in S :x^{-1}A \in q\} \in p$, where 
$x^{-1}A=\{ y \in S: x\cdot y \in A \}$. A non-empty subset $I$ of a semigroup $(T,\cdot)$ 
is called a left ideal of $S$ if $T\cdot I \subseteq I$, a right ideal of $S$ if 
$I\cdot T \subseteq I$ and a two sided ideal (or simply an ideal) if it is both a 
left and a right ideal. A maximal left ideal is a left ideal that does not 
contain any proper left ideal. Similarly we can define minimal right ideal 
and smallest ideal. Any compact Hausdorff right topological semigroup $(T,\cdot)$ 
has a smallest two sided ideal.\\
$$\begin{array}{ccc}
K(T) & = & \bigcup\{L:L \text{ is a minimal left ideal of } T\} \\
& = & \,\,\,\,\,\bigcup\{R:R \text{ is a minimal right ideal of } T\}\\
\end{array}$$
\\

We now present Bergelson's Characterizations of Central sets.\\

\textbf{Definition 1.2.} Let $S$ be a discrete semigroup and let $C$ be a subset of $S$. 
Then $C$ is centers if and only if there is an idempotent $p$ in $K(\beta S)$ such 
that $C \in p$.\\

To give a dynamical characterization of central set in arbitrary semigroup $S$, we 
need the following definitions from now on, $P_f(X)$ is the set of finite non-empty 
subsets of $X$, for any set $X$.\\

\textbf{Definition 1.3.}([Definition 3.1, 6]) Let $S$ be a semigroup and 
let $A\subseteq S$.

(a) The set $A$ is syndetic if and only if there is some $G \in P_f(S)$ such that 
    $S=\cup_{t \in G}t^{-1}A $.

(b) the set $A$ is piece-wise syndetic if and only if there is some  $G \in P_f(S)$ 
such that for an $F \in P_f(S)$ there is some $x \in S$ with 
$Fx\subseteq \cup_{t \in G}t^{-1}A $\\

Recall the definitions of proximality and uniform recurrence in a dynamical system 
from [3, Definition 1.2(b)] and [3, Definition 1.2(c)]\\

\textbf{Definition 1.4.} Let $(X,\langle T_s\rangle _{s \in S})$ is a dynamical system.

(a) A point $y \in S$ is uniformly recurrent if and only if , for every neighborhood 
    $U$ of $y$, $\{ s\in S : T_s(y) \in U  \}$ is syndetic.

(b) The points $x$ and $y$ of $X$ are proximal if and only if for every neighborhood 
    $U$ of the diagonal in $X\times X$, there is some $s \in S$ such that 
    $(T_s(x), T_s(y)) \in U$.\\
    
By [Theorem 2.4, 10], a subset $C$ of a semigroup $S$ is central if and only if 
there exist a dynamical system $(X,\langle T_s\rangle _{s \in S})$, points $x$ and $y$ of $X$ 
and a neighborhood $U$ of $y$ such that $y$ is uniformly recurrent, $x$ and $y$ 
are proximal and $C=\{s  \in S: T_s(x) \in U  \}$.\\

We now state basic definitions, conventions and results for dynamical 
characterization of members of certain idempotent ultrafilters.\\

\textbf{Definition 1.5.}([8, Definition 2.1])Let $S$ be a nonempty discrete space and 
$\mathcal{K}$ is a filter on $S$.

(a) $\bar{\mathcal{K}}= \{p \in \beta S : \mathcal{K} \subseteq p  \}$

(b) $L(\mathcal{K})=\{A\subseteq S: S\setminus A \not \in \mathcal{K}\}$\\

As is well known, the function $\mathcal{K} \rightarrow \bar{\mathcal{K}}$ is a bijection 
from the collection of all filters on $S$ onto the collection of all compact 
subspace of $\beta S$ [7, Theorem 3.20]. We also have the following important 
theorem relating the above two concepts.
\\

\textbf{Theorem 1.6.} Let $S$ be a nonempty discrete space and $\mathcal{K}$ a filter on $S$.

(a) $\bar{\mathcal{K}}=\{ p \in \beta S: A \in L(\mathcal{K})$ for all $A \in p$ \} 

(b) Let $\beta \subseteq L(\mathcal{K})$ be closed under finite intersections then 
there exists a $p \in \beta S$ with $\beta \subseteq p \subseteq L(\mathcal{K})$.\\
\textbf{Proof.} Both of these assertions follows from [7,Theorem 3.11]\\

\textbf{Definition 1.7.}(8, Definition 3.1) Let $(X,\langle T_s\rangle _{s \in S})$ be a 
dynamical system, $x$ and $y$ points in $X$, and $\mathcal{K}$ a filter on $S$. 
The pair $(x,y)$ is called jointly $\mathcal{K}$-recurrent if and only if for 
every neighborhood $U$ of $y$ we have $\{ s \in S: T_s(x) \in U$ and 
$T_s(y) \in U \} \in L(\mathcal{K})$.\\

Following theorem is [8, Theorem 3.3].\\

\textbf{Theorem 1.8.} Let $(S,.)$ be a semigroup, let $\mathcal{K}$ be a filter on 
$S$ such that $\bar{\mathcal{K}}$ is a compact subsemigroup of $\beta S$, and let 
$A \subseteq S$. Then $A$ is a member of an idempotent in $\bar{\mathcal{K}}$ if 
and only if there exists a dynamics system $(X,\langle T_s\rangle _{s \in S})$ with points 
$x$ and $y$ in $X$ and there exists a neighborhood $U$ of $y$ such that the 
pair $(x,y)$ is jointly $\mathcal{K}$-recurrence and $A=\{ s \in S: T_s(x) \in U \}$.\\

\textbf{Definition 1.9.}([6, Definition 1.2]) Let S be a discrete semigroup and let 
$C$ be a subset of $S$. Then $C$ is quasi-central if and only if there is an idempotent 
$p$ in Cl $K(\beta S)$ such that $C \in p$.\\

Now recall [8, Definition 4.1 and Definition 4.4]\\

\textbf{Definition 1.10.} Let $S$ be a semigroup.

(a) for each positive integer $m$ put $J_m=\{ (t_1,t_2,...,t_m)\in 
\mathbb{N}^m:t_1<t_2<...<t_m \}$

(b) given $m \in \mathbb{N}, a \in S^{m+1}, t \in J_m$, and $t \in \tau $, put 
$x(m,a,t,f)=\prod^m_{i=1}(a(i)f(t_i))a(m+1) $ where $\tau=\mathbb{N_S}$

(c) We call a subset $A\subseteq S$, a $C$-set if and only if there exists 
functions $m:P_f(\tau)\rightarrow \mathbb{N}$,$\alpha \in X^{S^{m(F)+1}}_{F \in P_f(\tau)}$, 
and $\tau \in X^{J_{m(F)}}_{F \in P_f(\tau)}$ such that the following two 
statements are satisfied:

(1) If $F,G \in P_f(\tau)$ and $F\subsetneq G$ then $\tau(F)(m(F))<\tau(G)(1)$.

(2) Whenever $m \in \mathbb{N}, G_1, G_2,...,G_m$ is a finite sequence in 
$P_f(\tau)$ with $G_1 \subsetneq G_2 \subsetneq  ...\subsetneq G_m$ and for 
each $i \in \{ 1,2,...,m \}, f_i \in G_i$ then we have\\
$\prod^m_{i=1}x(m(G_i), \alpha (G_i), \tau (G_i),f_i)\in A$

(d) We call a subset $A\subseteq S$, a $J$-set if and only if for every 
$F \in P_f(\tau)$, there exist $m \in \mathbb{N},a \in S^{m+1}$ and 
$t \in J_m$ such that for all $f \in F$, $x(m,a,t,f)\in A$\\

\textbf{Definition 1.11.} Let $(X,\langle T_s\rangle _{s \in S})$ be a dynamical system and 
let $x,y \in X$.

(a) The pair $(x,y)$ is jointly intermittently uniformly recurrent 
(abbreviated JIUR) if and only if for every neighborhood $U$ of $y$, 
$\{ s \in S: T_s(x)\in U$ and $T_s(y)\in U \}$ is piecewise syndetic.

(b) the pair $(x,y)$ is jointly intermittently almost uniform recurrent 
(abbreviated as JIAUR) if and only if for every neighborhood $U$ of $y$, 
$\{ s \in S: T_s(x)\in U$ and $T_s(y)\in U \}$ is a $J$-set.\\

Using theorem 1.8, we have dynamical characterizations of quasi-central set 
and $C$-set in terms of JIUR and JIAUR respectively.\\

\textbf{Theorem 1.12.} Let S be a semigroup and let $C\subseteq S$. The set $C$ 
is quasi-central if and only if there exists a dynamical system $(X,\langle T_s\rangle _{s \in S})$, 
points $x$ and $y$ of $X$ such that $x$ and $y$ are JIUR, and a neighborhood $U$ 
of $y$ such that $C=\{ s \in S: T_s(x)\in U\}$.\\

\textbf{Theorem 1.13.} Let S be a semigroup and let $C\subseteq S$. The set $C$ 
is $C$-set if and only if there exists a dynamical system $(X,\langle T_s\rangle _{s \in S})$, 
points $x$ and $y$ of $X$ such that $x$ and $y$ are JIAUR, and a neighborhood $U$ 
of $y$ such that $C=\{ s \in S: T_s(x)\in U  \}$.\\

Now we will be considering semigroups which are dense in $((0,\infty),+)$. 
Here 'dense' means with respect to the usual topology on $((0,\infty),+)$.\\

\textbf{Definition: 1.14.} If $S$ be a dense subsemigroup of $((0,\infty),+)$ 
then $O^{+}(S)=\{ p \in \beta S :$ for an $ \epsilon >0, S\cup(0,\epsilon) \in p\}$\\

It was proved in [5, Lemma 2.5] that $O^{+}(S)$ is a compact right topological 
subsemigroup of $(\beta S,+)$. It was also noted there in $O^{+}(S)$ is disjoint 
from $K (\beta S)$ and hence gives some new information which is not available 
from $K (\beta S)$.\\

Being a compact right topological semigroup , $O^{+}(S)$ has a minimal ideal 
$K (O^{+}(S)$.\\

Like discrete case we can define central set, quasi-central set, $C$-set near zero. 
Dynamical characterization of central set near zero, quasi central set near zero 
and $C$-set near zero is established in section 2, section 3 and section 4 respectively..\\

\section {dynamical characterization of central set near zero}

Let us start this section with the following well-known definition of central 
set near zero[5, Definition 4.1(a)]\\

\textbf{Definition 2.1.} Let $S$ be a dense subsemigroup of $((0,\infty),+)$. 
A set $C\subseteq S$ is a central near zero if and only if there is some 
idempotent $p\in K(O^{+}(S))$ with $C \in p$\\

Following definition is [5, Definition 3.2(b)].\\

\textbf{Definition 2.2.} Let $S$ be a dense subsemigroup of $((0,\infty),+)$. A subset 
$B$ of $S$ is syndetic near zero if and only if for every $\epsilon>0$, there exists 
some $F \in P_f((0,\epsilon)\cap S)$ and some $\delta>0$ such that 
$S\cap (0,\delta)\subseteq \cup_{t \in F}(-t+B)$.\\

We shall now introduce the the notion of uniform recurrence and proximality near zero.\\

\textbf{Definition 2.3.} Let S be a dense subsemigroup of $((0,\infty),+)$ and 
$(X,\langle T_s\rangle _{s\in S})$ be a topological dynamical system.

(a) A point $x \in X$ is a uniformly recurrent point near zero if and only if 
for each neighborhood $W$ of $x$, $\{s \in S: T_s \in W  \}$ is syndetic near zero.

(b) Points $x$ and $y$ of $X$ are proximal near zero if and only if for every 
neighborhood $U$ of the diagonal in $X\times X$, for each $\epsilon>0$ there 
exists $s\in S \cap (0, \epsilon)$ such that $(T_s(x), T_s(y)) \in U$.\\

We now recall [7, Theorem 19.11].\\

\textbf{Theorem 2.4.} Let $(X,\langle T_s \rangle_{s\in S})$ be a dynamical system and define 
$\theta : S \rightarrow X^X$ by $\theta(s)=T_s$. Then $\tilde{\theta}$ is a 
continuous homomorphism from $\beta S$ onto the enveloping semigroup of 
$(X,\langle T_s\rangle _{s\in S})$. ($\tilde{\theta}$ is the continuous extension of $\theta$)\\

The following notation will be convenient in the next section.\\

\textbf{Definition 2.5.}([7, Definition 19.12]) Let $(X,\langle T_s\rangle _{s\in S})$ be a dynamical 
system and define $\theta : S\rightarrow X^X$ by $\theta(s)=T_s$. For each 
$p \in \beta S$, let $T_p= \tilde{\theta}(p)$.\\

As a immediate consequence of theorem 2.4 we have the following remark [7, Remark 19.13]\\

\textbf{Remark 2.6.} Let $(X,\langle T_s\rangle _{s\in S})$ be a dynamical system and let 
$p,q \in \beta S$. Then $T_p\circ T_q=T_{pq}$ and for each 
$x \in X$, $T_p(x)=p-lim_{s \in S}T_s(x)$.\\

Clearly it is easy to see that, points $x$ and $y$ of $X$ are proximal near zero 
if and only if there is some $p \in O^{+}(S)$ such that $T_p(x)=T_p(y)$ \\

\textbf{Lemma 2.7.} Let S be a dense subsemigroup of $((0,\infty),+)$. Let $(X,\langle T_s\rangle _{s\in S})$ be a topological dynamical system and $L$ be a minimal left ideal of $O^{+}(S)$ and $x \in X$.

The following statements are equivalent :

(a) The point $x$ is a uniformly recurrent point near zero of $(X,\langle T_s\rangle _{s\in S})$.

(b) There exists $u \in L$ such that $T_u(x)=x$.

(c) There exists $y \in X$ and an idempotent $u \in L$ such that $T_u(y)=x$.

(d)there exists an idempotent $u \in L$ such that $T_u(x)=x$.\\

\textbf{Proof.} (a)$\Rightarrow$(b)
Choose any $v \in L$. Let $N$ be a set of neighborhoods of $x$ in $X$. 
For each $U \in N$ let $B_U=\{ s \in S:T_s(x) \in U  \}$.
Since $x$ is uniformly recurrent point near zero, each $B_U$ is syndetic near zero, 
for every $\epsilon >0$ there is some $F_{\epsilon} \in P_f((0,\epsilon)\cap S)$ and 
some $\delta>0$ such that $S\cup(0,\delta)\subseteq \cup_t \in F_{\epsilon}(-t+B)$. 
So, for each $U \in N$ and $\epsilon>0$, pick $t_{(U,\epsilon)}\in F_{(U,\epsilon)}$ 
such that

$-t_{(U,\epsilon)}+B_U\in v$
Given $U \in N$ and $\epsilon > 0$, let 
 
$C_{(U,\epsilon)}=\{ t_{(v,\epsilon)}: v\in N$ and $V\subseteq U  \}$ and 
$C_U= \cup_{\epsilon >0}C_{(U,\epsilon)}$, then

$\{ C_U:U \in N  \} \cup \{ S \cap (0,\epsilon): \epsilon>0  \}$ has the finite 
intersection property.
 
Now pick $w \in O^{+}(S)$ such that $\{C_U: U \in N   \}\subseteq w$ and let 
$u=w+v$. Then $u \in L$ since $L$ is a left ideal of $O^{+}(S)$.
 
To see that $T_U(x)=x$, let $U \in N$. We need to show that $B_U \in u $, suppose 
instead that $B_U \not\in u$. 

Then $\{ t \in S: -t+B_U \not\in v  \}$ and $C_U \in w$ and so pick $t \in C_U$ 
such that, $-t+B_U\not\in v$. Pick $V\in N$ with $V\subseteq U$ such that 
$t=t_{(V,\epsilon)}$ for some $\epsilon >0$. Then $-t+B_v \in v$ and $-t+B_V \subseteq -t+B_v$, a contradiction.\\
(b)$ \Rightarrow (c)$ Let $K=\{ v \in L:T_v(x)=x  \}$. It suffices to show that $K$ is 
a compact subsemigroup of $L$, since then $K$ has an idempotent. By assumption, $K\neq \phi$ . 
Further if $v \in L \setminus K$, then there is some neighborhood $U$ of $x$ such that 
$B=\{ s \in S:T_s(x)\in U   \}\not\in v$. Then Cl$B$ is a neighborhood of $v$ in $\beta S$ 
which misses $K$. Finally, to see that $K$ is a semigroup, let $v,w\in K$. Then by Remark 2.6, 
$T_{v+w}(x)=T_v(T_w(x))=T_v(x)=x$.\\
(c)$\Rightarrow$(d)\\
Again we use remark 2.6: $T_u(x)=T_u(T_u(y))=T_{u+u}(y)=T_u(y)=x$\\
(d)$\Rightarrow$(a)\\
Let $U$ be a neighborhood of $x$ and let $B=\{ s \in S:T_s(x)\in U\}$ and suppose 
that $B$ is not syndetic near zero. Then there exists $\epsilon >0$ such that 
$\{ S\setminus \cup_{t \in F}(-t+B)$: F is a finite nonempty subset of 
$S \cap (0,\epsilon)  \}$ $\cup$ $ \{S \cap (0,\delta): 0<\delta < \infty\}  $
has the finite intersection property. So pick some $w \in O^{+}(S)$ such that 
$\{ S\setminus \cup_{t \in F}(-t+B)$: F is a finite nonempty subset of 
$S \cup (0,\epsilon)\}\subseteq w$.\\
Then $(O^{+}(S)+w)\cap ClB=\phi $(For suppose instead one had some 
$v \in O^{+}(S)$ with $B \in v+w$. Then pick some $t \in P_f(S \cap (0, \epsilon)$ 
with $-t+B \in w$  ).\\
Let $L^{'}= O^{+}(S)+w$. Then $L^{'}$ is a left ideal of $O^{+}(S)$, so $L^{'}+u$ is a 
left ideal of $O^{+}(S)$ which is contained in $L$ , and hence $L^{'}+u=L$. Thus we 
may pick some $v \in L^{'}$ such that $v+u=u$. Again using Remark 2.6, 
$T_v(x)=T_v(T_u(x))=T_{v+u}(x)=T_u(x)=x$, so in particular $B \in v$. But, $v \in L^{'}$ 
and $L^{'}\cap ClB=\phi$, a contradiction.\\

\textbf{Lemma 2.8.} Let $S$ be a dense a subsemigroup of $((0,\infty),+)$. 
Let, $(X,\langle T_s\rangle _{s \in S})$ be a topological dynamical system and let $x \in X$. 
Then for each $\epsilon > 0$ there is a uniformly recurrent point near zero 
$y \in Cl\{ T_s(x): s \in S \cap (0, \infty)   \}$ such that $x$ and $y$ are 
proximal near zero.\\
\textbf{Proof.} Let $L$ be any minimal left ideal of $O^{+}(S)$ and pick an 
idempotent $u \in L$. Let $y=T_u(x)$. For each $\epsilon > 0$, clearly  
$y \in Cl\{ T_s(x): s \in S \cap (0, \infty)   \}$. By Lemma 2.7, $y$ is a 
uniformly recurrent point near zero of $(X,\langle T_s\rangle _{s \in S})$. By Remark 2.6 we have $T_u(y)=T_u(T_u(x))=T_{u+u}(x)=T_u(x)$ . So $x$ and $y$ are proximal near zero.\\

\textbf{Lemma 2.9.} Let $S$ be a dense subsemigroup of $((0, \infty),+)$. 
Let $(X,\langle T_s\rangle _{s \in S})$ be a topological dynamical system and let $x,y \in X$. 
If $x$ and $y$ are proximal near zero, then there is a minimal left ideal 
$L$ of $O^{+}(S)$ such that $T_u(x)=T_u(y)$ for all $u \in L$\\
\textbf{Proof.} Pick $v \in O^{+}(S)$ such that $T_v(x)=T_v(y)$ and pick a 
minimal left ideal $L$ of $O^{+}(S)$ such that $L \subseteq O^{+}(S)+v$. to 
see that $L$ is as required , let $u \in L$ and choose $w \in O^{+}(S)$ such 
that $u=w+v$. Then again using Remark 2.6, we have 
$T_u(x)=T_{w+v}(x)=T_w(T_v(x)=T_w(T_v(y))=T_{w+v}(y)=T_u(y)$\\

\textbf{Lemma 2.10.} Let $S$ be a dense subsemigroup. Let $(X,\langle T_s\rangle {s \in S})$ be 
a topological dynamical system and let $x,y \in X$. There is an idempotent $u$ in 
$K(O^(S))$ such that $T_u(x)=y$ if and only if both $y$ is uniformly recurrent 
near zero and $x$ and $y$ are proximal near zero.\\
\textbf{Proof.} 
$(\Rightarrow)$.
Since $u$ is a minimal idempotent of $O^{+}(S)$, there is a minimal left ideal 
$L$ of $O^{+}(S)$ such that $u \in L$. Thus by Lemma 2.7, $y$ is uniformly 
recurrent near zero. By Remark 2.6, 
$T_u(y)=T_u(T_u(x))=T_{u+u}(x)=T_u(x)$. so $x$ and $y$ are proximal near zero.\\
$(\Leftarrow)$ Pick by Lemma 2.9 a minimal ideal $L$ of $O^{+}(S)$ such that 
$T_u(x)=T_u(y)$ for all $u \in L$. Pick by Lemma 2.7 an idempotent $u \in L$ 
such that $T_u(y)=y$.\\

We now give a dynamical Characterization of Central sets near zero in the 
following theorem.\\

\textbf{Theorem 2.11.} Let $S$ be a dense subsemigroup of $((0, \infty),+)$ and let 
$B \subseteq S$. Then $B$ is central near zero if and only if there exists a 
topological dynamical system $(X,\langle T_s\rangle {s \in S})$ and there exists $x,y \in X$ 
and a neighborhood $U$ of $y$ such that $x$ and $y$ are proximal near zero , $y$ 
is uniformly recurrent near zero, and $B=\{ s \in S: T_s(x) \in U   \}$.\\
\textbf{Proof.} $(\Rightarrow)$ Let $G=S\cup \{0\}$ and $X=\prod_{s \in G}\{0,1\}$ and 
for $s \in S$ define $T_s:X\rightarrow X$ by $T_s(x)(t)=x(x+t)$ for all $t \in G$. It 
is easy to see that $T_s$ is continuous. Now let $x=\chi _{B}$, the characteristic 
function of $B$. That is, $x(t)=1$ if and only if $t \in B$. Pick a minimal idempotent 
in $O^{+}(S)$ such that $B \in u$ and let $y=T_u(x)$. Then by Lemma 2.10, $y$ is 
uniformly recurrent near zero and $x$ and $y$ are proximal near zero.\\
Now let $U=\{ x \in X:z(0)=y(0) \}$. Then $U$ is a neighborhood of $y$ in $X$. 
We note that $y(0)=1$. Indeed, $y=T_u(x)$ so, $\{ s \in S:T_s(x) \in U  \}\in u$ 
and we may choose some $s \in B$ such that $T_s(x) \in U$. Then $y(0)=T_s(x)(0)=x(s+0)=1$
Thus given any $s \in S$, 
\\
$s \in B\Leftrightarrow x(s)=1 \Leftrightarrow T_s(x)(0)=1 \Leftrightarrow T_s(x)\in v$
\\
$(\Leftarrow)$ Choose a topological dynamical system $(X,\langle T_s\rangle {s \in S})$, points $x,y \in X$ and a neighborhood $U$ of $y$ such that $x$ and $y$ are proximal, $y$ is uniformly recurrent and $B= \{ s \in S : T_s(x) \in U  \}$. Choose by Lemma 2.10 a minimal idempotent $u$ in $O^{+}(S)$ such that $T_u(x)=y$. Then $B \in u$\\

\section{Dynamical characterization of quasi-central set near zero}
\textbf{Definition 3.1.} Let $S$ be a dense subsemigroup of $((0, \infty), +)$. 
Then $C$ is said to be quasi central near zero if and only if there is an idempotent 
$p$ in $Cl\text{ }K(O^+(S))$ such that $C \in p$.\\ 

Quasi-central sets near zero have some significant virtues. In the first place it satisfies conclusion of central sets theorem near zero. Secondly, in [5] and [9] combinatorial characterization of central set near zero and quasi-central set near zero are obtained respectively. The characterization of quasi-central set near zero is much simpler than the characterization of Central sets near zero.\\

In this section we shall deduce the dynamical characterization of quasi-central set near zero.\\

For this we need the following two definitions.\\

\textbf{Definition 3.2.}([5, Definition $3.4$])
Let $S$ be a dense subsemigroup of $((0,\infty), +)$. A subset $A$ of $S$ is
piecewise syndetic near zero if and only if there exists sequences 
$\langle F_n\rangle_{n=1}^{\infty}$ and $\langle \delta_n\rangle_{n=1}^{\infty}$ such that\\ 
(1) for each $n \in \mathbb{N}$, $F_n \in \mathcal{P}_f((0,\frac{1}{n}) \cap S)$ 
and $\delta_n \in (0, \frac{1}{n})$ and\\
(2) for all $G \in \mathcal{P}_f(S)$ and $\mu > 0$ there is some $x \in (0,\mu) \cap S$ 
such that for all $n \in \mathbb{N}$, we have  
$$(G \cap (0,\delta_n)) + x \subseteq \cup_{t \in F_n}(t+A).$$\\

As in the discrete case let us now introduce the notion of jointly 
interminittently uniform recurrent near zero.\\

\textbf{Definition 3.3.} Let $(X, \langle T_s\rangle_{s\in S})$ be a dynamical system 
and let $x,y \in X$. The pair $(x,y)$ is jointly interminittently uniformly recurrent 
near zero (abbriviated as $JIUR_0$) if and only if for every neighbourhood $U$ of $y$, 
the set $\{s \in S : T_s(x) \in U \text{ and } T_s(y)\in U\}$ is piecewise syndetic near zero.\\

For our purpose, we state the following theorem[5, Theorem $3.5$]\\

\textbf{Theorem 3.4.} Let $S$ be a dense subsemigroup of $((0,\infty), +)$ and let $A \subseteq S$. 
Then $K(O^+(S)) \cap Cl \text{ }A \neq \phi$ if and only if $A$ is piecewise syndetic near zero.\\

\textbf{Lemma 3.5.}
Let $S$ be a dense subsemigroup of $((0, \infty), +)$ and $$\mathcal{K} = \{A \subseteq S : S \setminus A \text{ is not piecewise syndetic near zero}\}.$$ Then $\mathcal{K}$ is a filter on $S$ with $Cl \text{ } K(O^+(S)) = \overline{\mathcal{K}}$, which is a compact subsemigroup of $\beta S$.\\
\textbf{Proof.} It is a routine exercise to show that $\mathcal{K}$ is non-empty, does not contain the empty set and is closed under super set. To show that $\mathcal{K}$ is a filter it is enough to prove that $\mathcal{K}$ is closed under finite intersection. Let $A, B \in \mathcal{K}$, then both $S \setminus A$ and $S \setminus B$ are piecewise syndetic near zero. Now we shall show that $A \cap B \in \mathcal{K}$ i.e. $S \setminus (A \cap B)$ is not piecewise syndetic near zero. If possible let $S \setminus (A \cap B)$ is piecewise syndetic near zero. So by Theorem $3.4$, there exists $p \in Cl \text{ }K(O^+(S))$ such that $S \setminus (A \cap B) \in p$. This $(S \setminus A) \cup (S \setminus B) \in p$ and therefore $S \setminus A \in p$ or $S \setminus B \in p$, which is a contradiction. So our claim is proved.

Observe that under the assumption that $\mathcal{K}$ is a filter, $$\mathcal{L}(\mathcal{K}) = \{A \subseteq S : A \text{ is piecewise syndetic near zero}\}$$
So by Theorem 1.6 and Theorem 3.4, we have $Cl \text{ }K(O^+(S)) = \overline{\mathcal{K}}$, which is a compact subsemigroup of $\beta S$.\\

In the following theorem we shall give the dynamical characterization of quasi-central set near zero.\\
\textbf{Theorem 3.6.} Let $S$ be a dense subsemigroup of $((0, \infty), +)$ and let $A \subseteq S$. The set $A$ is quasi-central near zero if and only if there exists a dynamical system $(X,\langle T_s \rangle_{s \in S})$, points $x$ and $y$ of $X$ such that $x, y$ are $JIUR_0$ and a neighbourhood $U$ of $y$ such that $$A = \{s \in S : T_s(x) \in U\}.$$\\
\textbf{Proof.} Let $$\mathcal{K} = \{B \subseteq S : B \text{ is not a piecewise syndetic set near zero}\},$$ and note that $$\mathcal{L}(\mathcal{K}) = \{A \subseteq S : A \text{ is piecewise syndetic near zero}\}.$$ By Lemma 3.5, we have $\mathcal{K}$ is a filter and $\overline{\mathcal{K}} = Cl \text{ } K(O^+(S))$ which is a compact subsemigroup of $\beta S$. Now choose an idempotent $p$ in $\overline{\mathcal{K}} = Cl \text{ } K(O^+(S))$ such that $c \in p$. Now we can apply Theorem $1.8$ to prove our required statement.\\

\section{Dynamical characterization of $C$-set near zero}
We start by giving the combinatorial definitions of $C$-set near zero. As this combinatorial definition is rather complicated, we shall soon state an algebraic characterization showing that $C$-sets are members of idempotents in a certain compact subsemigroup.\\ 
\textbf{Definition 4.1}([1, Definition 3.1])
Let $S$ be a dense subsemigroup of $((0,\infty), +)$. The set of sequences in $S$ converging to $0$ is denoted by $\tau_0$.\\

\textbf{Definition 4.2.}([1, Definition $3.6(a)$])
Let $S$ be a dense subsemigroup of $((0,\infty),+)$ and let $A \subseteq S$. We say $A$ is a $C$-set near zero if and only if for each $\delta \in (0,1)$, there exists functions $a_{\delta} : \mathcal{P}_f(\tau_0) \to \mathcal{P}_f(\mathbb{N})$, such that\\
(1) $\alpha_{\delta}(F) < \delta$ for each $F \in \mathcal{P}_f(\tau_0)$,
(2) if $F, G \in \mathcal{P}_f(\tau_0)$ and $F \leq G$, then $\max H_{\delta}(F) < \min H_{\delta}(G)$ and 
(3) whenever $m \in \mathbb{N}$, $G_1,G_2, \cdots , G_m \in \mathcal{P}_f(\tau_0)$, $G_1 \subseteq G_2\subseteq \cdots \subseteq G_m$ and for each $i \in \{1,2, \cdots ,m\}$, $f_i \in G_i$, one has $$\sum_{i=1}^{m}\big(\alpha_{\delta}(G_i)+ \sum_{t \in H_{\delta}(G_i)}f_i(t)\big) \in A.$$\\

We now recall [1, Definition $3.2$]  and [1, Definition $3.6(b)$].\\

\textbf{Definition 4.3.}
Let $S$ be a dense subsemigroup of $((0, \infty),+)$ and let $A \subseteq S$.
(1) $A$ is said to be $J$-set near zero if and only if whenever $F \in \mathcal{P}_f(\tau_0)$ and $\delta > 0$, there exists $a S \cap (0, \delta)$ and $H \in \mathcal{P}_f(\mathbb{N})$ such that for each $f \in F$, $a + \sum_{t \in H}f(t) \in A$.
(2) $J_0(S) = \{p \in O^+(S) : \text{ for all } A \in p, \text{ is a J-set near zero}\}$.\\

\textbf{Lemma 4.4.}
Let $S$ be a dense subsemigroup of $((0, \infty),+)$ and $A_1$, $A_2$ are subsets of $S$. If $A_1 \cup A_2$ is $J$-set near zero then either $A_1$ or $A_2$ is a $J$-set near zero.\\
\textbf{proof} See [1, Lemma $3.8$].\\

\textbf{Definition 4.5.} Let $(X, \langle T_s\rangle_{s \in S})$ be a dynamical system and $x,y \in X$. The pair $(x,y)$ is jointly almost uniform recurrent (abbreviated $JIAUR_0$) if and only if for every neighbourhood $U$ of $y$, $\{s \in S : T_s(x) \in U \text{ and } T_s(y) \in U\}$ is a $J$-set near zero.\\

\textbf{Lemma 4.6.}
Let $S$ be a dense subsemigroup of $((0,\infty),+)$ and $$\mathcal{K}= \{A \subseteq S : S \setminus A \text{ is not a }J\text{-set near zero}\}.$$ Then $\mathcal{K}$ is a filter on $S$ with $J_0(S) = \overline{\mathcal{K}}$ and $J_0(S)$ is a compact subsemigroup of $\beta S$.\\
\textbf{Proof.} It is easy to see that $\mathcal{K}$ is non-empty, does not contain the empty set and is closed under super sets. By Lemma 4.4, it follows that $\mathcal{K}$ is closed under finite intersection.

Under the assumption that $\mathcal{K}$ is a filter, we have $$\mathcal{L}(\mathcal{K}) = \{ A \subseteq S : A \text{ is a }J\text{-set}\}.$$ From Theorem $1.6$ it follows that $J_0(S) = \overline{\mathcal{K}}$. Finally, the fact that $J_0(S)$ is a subsemigroup of $\beta S$ follows from [1, Theorem $3.9$].\\

\textbf{Lemma 4.7.} Let $S$ be a dense subsemigroup of $((0,\infty),+)$ and $A \subseteq S$. Then $A$ is a $C$-set near zero if and only if there exists an idempotent $p \in J_0(S)$ such that $A \in p$.\\
\textbf{Proof} This is proved in [1, Theorem $3.14$].\\

The following theorem gives a dynamical characterization of $C$-set near zero.\\

\textbf{Theorem4.8} Let $S$ be a dense subsemigroup of $((0,\infty),+)$ and $A \subseteq S$. Then $A$ is a $C$-set near zero if and only if there exists a dynamical system $(X, \langle T_s \rangle_{s \in S})$ with points $x$ and $y$ in $X$ such that $x$, $y$ are $JIAUR_0$ and a neighbourhood $U$ of $y$ such that $A = \{s \in S : T_s(x) \in U\}$.\\
\textbf{Proof} Let $\mathcal{K} = \{B \subseteq S : S \setminus B \text{ is not a }J \text{-set near zero}\}$. Since Lemma 4.7, characterizes $C$-set near zero in terms of idempotents in $\overline{\mathcal{K}}$, we can apply Theorem  $1.8$ to prove our statement.\\

\textbf {Acknowledgement.}  The author is grateful to Prof.
Swapan Kumar Ghosh of Ramakrishna Mission Vidyamandira for 
continuous inspiration and a number of valuable suggestions 
towards the improvement of the paper.

\end{document}